\documentclass[12pt]{amsart}
\usepackage{amsmath, amssymb, latexsym}
\usepackage{mathrsfs}
\usepackage{enumerate}
\usepackage{color}
\usepackage{graphics}
\usepackage{verbatim}
\usepackage[colorlinks=true, pdfstartview=FitV, linkcolor=blue,%
citecolor=blue, urlcolor=blue]{hyperref}
\usepackage[all, knot]{xy}
\xyoption{arc}

\newtheorem{Thm}{Theorem}[section]
\newtheorem{Proposition}[Thm]{Proposition}
\newtheorem{Corr}[Thm]{Corollary}
\newtheorem{Lem}[Thm]{Lemma}
\newtheorem{Defn}[Thm]{Definition}

\theoremstyle{definition}

\newtheorem{Exam}[Thm]{Example}
\newtheorem{Remark}[Thm]{Remark}

\newcommand{\nc}{\newcommand}

\nc{\Lemma}{\begin{Lem}}
\nc{\enlemma}{\end{Lem}}
\nc{\Cor}{\begin{Corr}}
\nc{\encor}{\end{Corr}}
\nc{\Th}{\begin{Thm}}
\nc{\entheorem}{\end{Thm}}
\nc{\Prop}{\begin{Proposition}}
\nc{\enprop}{\end{Proposition}}
\nc{\Def}{\begin{Defn}}
\nc{\enDef}{\end{Defn}}
\nc{\Ex}{\begin{Exam}}
\nc{\enEx}{\end{Exam}}
\nc{\Rem}{\begin{Remark}}
\nc{\enRem}{\end{Remark}}

\numberwithin{equation}{section}

\renewcommand{\to}[1][]{\xrightarrow{#1}{}}

\newcommand{\Z}{\mathbf{Z}}
\newcommand{\Q}{\mathbf{Q}}
\newcommand{\R}{\mathcal{R}}

\newcommand{\A}{\mathbf{A}}
\nc{\sho}{\mathscr{O}}

\newcommand{\g}{\mathfrak{g}}

\newcommand{\wt}{{\rm wt}}

\nc{\on}{\operatorname}

\newcommand{\Hom}{\on{Hom}}

\newcommand{\ch}{{\rm ch}}

\newcommand{\End} { {\rm End}}
\newcommand{\id} { {\rm id}}

\nc{\cor}{\mathbf{k}}
\nc{\KLR}{quiver Hecke algebra}
\nc{\KLRs}{quiver Hecke algebras}
\nc{\seteq}{\mathbin{:=}}
\newcommand{\soplus}{\mathop{\mbox{\normalsize$\bigoplus$}}\limits}
\newcommand{\ssum}{\mathop{\mbox{\footnotesize$\displaystyle\sum$}}\limits}
\nc{\cl}{\colon}
\nc{\set}[2]{\left\{#1\mid #2\right\}}
\nc{\Id}{\operatorname{Id}}
\nc{\Ker}{\on{Ker}}
\nc{\Coker}{\on{Coker}}
\nc{\coh}{\mathrm{coh}}
\nc{\Mod}{\on{Mod}}
\nc{\Modc}{\on{Mod_\coh}}
\nc{\Proj}{\on{Proj}}
\nc{\Rep}{\on{Rep}}
\newcommand{\isoto}[1][]{\mathop{\xrightarrow[#1]%
{{\raisebox{-.6ex}[0ex][-.6ex]{$\mspace{2mu}\sim\mspace{2mu}$}}}}}

\nc{\To}[1][\quad]{\to[\;#1\;]}
\nc{\hs}{\hspace*}
\nc{\vs}{\vspace*}
\nc{\bF}{\overline{F}}
\nc{\epi}{\twoheadrightarrow}
\nc{\mono}{\rightarrowtail}
\nc{\be}{\begin{enumerate}}
\nc{\ee}{\end{enumerate}}
\nc{\ba}{\begin{array}}
\nc{\ea}{\end{array}}
\nc{\eq}{\begin{eqnarray}}
\nc{\eneq}{\end{eqnarray}}
\nc{\eqn}{\begin{eqnarray*}}
\nc{\eneqn}{\end{eqnarray*}}
\nc{\ran}{\rangle}
\nc{\lan}{\langle}
\nc{\bl}{\bigl(}
\nc{\br}{\bigr)}
\nc{\bnum}{\be[{\rm(i)}]}
\nc{\enum}{\ee}
\nc{\bna}{\be[{\rm(a)}]}
\nc{\Proof}{\begin{proof}}
\nc{\QED}{\end{proof}}
\newcommand{\scbul}{{\,\raise1pt\hbox{$\scriptscriptstyle\bullet$}\,}}
\nc{\tens}{\mathop\otimes}
\nc{\E}[1][i]{{\mathsf{E}_{#1}^\Lambda}}
\nc{\F}[1][i]{{\mathsf{F}_{#1}^\Lambda}}
\nc{\x}[1][{1}]{a^\Lambda(x_{#1})}
\nc{\vphi}{\varphi}
\nc{\haut}{\mathrm{ht}}
\nc{\al}{\alpha}
\nc{\La}{\Lambda}
\nc{\Gr}{\on{Gr}}
\nc{\la}{\lambda}
\nc{\noi}{\noindent}
\nc{\eps}{\varepsilon}
\nc{\RL}{R^\La}
\nc{\One}{\mathbf{1}}
\nc{\comout}{}
\nc{\bigwr}{\mbox{\large$\wr$}}
\nc{\heps}{\widehat{\eps}}
\nc{\heta}{\widehat{\eta}}
\nc{\tE}{\widetilde{E}}
\nc{\tF}{\widetilde{F}}
\nc{\Fil}{\Gamma}
\nc{\op}{\mathrm{opp}}
\nc{\shc}{\mathscr{C}}
\nc{\Fct}{\on{Fct}}
\nc{\HH}{\mathsf{H}}
\nc{\nn}{\nonumber}
\nc{\gmod}{\mbox{-$\mathrm{gmod}$}}
\nc{\proj}{\mbox{-$\mathrm{proj}$}}
\nc{\cg}{{\mathrm{c}_{\mathrm{gen}}}}
\nc{\h}{\mathfrak{h}}
\nc{\rootl}{\mathsf{Q}}
\nc{\tf}{\tilde{f}}
\nc{\te}{\tilde{e}}
\nc{\cora}{{\cor(A)}}
\nc{\Par}{{X(A)}}
\nc{\Dual}{\mathrm{D}}
\nc{\ol}{\overline}
\nc{\Irr}{\on{\ol{Irr}}}
\nc{\CO}{\mathcal{O}}
\nc{\Res}{\on{Red}}
\nc{\G}{\mathrm{G}^{\mathrm{up}}}
\renewcommand{\L}{\mathcal{L}}
\nc{\Ln}{L}
\nc{\K}[1]{\mathrm{K}\bl #1\br }
\nc{\M}{\mathcal{M}}
\nc{\vac}{\phi}
\renewcommand{\hom}[1][]%
{{\mathscr{H}\mspace{-4mu}om}_{\raise1.5ex\hbox to.1em{}#1}}
\newcommand{\hhom}[1][]%
{{\mathscr{H}\mspace{-4mu}om}_{\raise1.5ex\hbox to.1em{}#1}}

\newlength{\my}
\begin{document}

\title[Notes on parameters of quiver Hecke algebras]
{Notes on parameters of quiver Hecke algebras}
\author[Masaki Kashiwara]{Masaki Kashiwara}
\thanks{This work was supported by Grant-in-Aid for
Scientific Research (B) 22340005, Japan Society for the Promotion of Science.}
\address{Research Institute for Mathematical Sciences, Kyoto University, Kyoto 606-8502, Japan}
\email{masaki@kurims.kyoto-u.ac.jp}

\date{April 7, 2012}


\begin{abstract}
Varagnolo-Vasserot and Rouquier proved that,
 in a symmetric generalized Cartan matrix case,
the simple modules over
the quiver Hecke algebra
with a special parameter
correspond to the upper 
global basis.
In this note we show that the simple modules over
the quiver Hecke algebras
with a {\em generic} parameter also
correspond to the upper global basis
in a symmetric generalized  Cartan matrix case.
\end{abstract}

\maketitle


\section{Introduction}

Lascoux-Leclerc-Thibon (\cite{LLT}) conjectured that
the irreducible representations of Hecke algebras of type $A$ are
controlled by the upper global basis (\cite{Kash91,Kash93})
(or dual canonical basis (\cite{Lus93})) of
the basic representation of the affine quantum group $U_q(A^{(1)}_\ell)$.
Then Ariki (\cite{A}) proved this conjecture by generalizing it to 
cyclotomic Hecke algebras.
The crucial ingredient in his proof
was the fact that the cyclotomic Hecke algebras
categorify the irreducible highest weight representations
of $U(A^{(1)}_\ell)$. Because of the lack of grading on
the cyclotomic Hecke algebras, these algebras do not
categorify the representation of the quantum group.

Then Khovanov-Lauda and Rouquier 
introduced independently
a new family of {\em graded} algebras, a generalization of
affine Hecke algebras of type $A$,
in order to categorify arbitrary quantum groups (\cite{KL09, KL08, R08}).
These algebras are  called  {\em Khovanov-Lauda-Rouquier algebras} or 
{\em quiver Hecke algebras.}

 Let $U_q(\g)$ be the quantum
group associated with a symmetrizable Cartan datum and let
$\{R(\beta)\}_{\beta \in \rootl^{+}} $ be the corresponding
\KLRs. Then it was shown in \cite{KL09,
KL08} that there exists an algebra isomorphism
$$U_{\A}^{-}(\g) \simeq \bigoplus_{\beta \in
\rootl^{+}}\K{R(\beta)\proj},$$ where $U_{\A}^-(\g)$ is the integral form
of the half $U_q^-(\g)$ of the quantum group 
$U_{q}(\g)$ with $\A = \Z[q, q^{-1}]$,
and $\K{R(\beta)\proj}$ is the Grothendieck group of 
the category $R(\beta)\proj$ of finitely generated
projective graded $R(\beta)$-modules. The positive 
root lattice is denoted by $\rootl^+$.
By the duality, we have
\eq U_{\A}^{-}(\g)^* \simeq \bigoplus_{\beta \in
\rootl^{+}}\K{R(\beta)\gmod},\label{eq:gmod}\eneq
where $U_{\A}^-(\g)^*$ 
is the direct sum of the dual of the weight space
$U_{\A}^{-}(\g)_{-\beta}$ of  $U_\A^-(\g)$, and 
$R(\beta)\gmod$ is the abelian category of graded $R(\beta)$-modules
which are finite-dimensional over the base field $\cor$.

When the generalized Cartan
matrix is a symmetric matrix, Varagnolo and Vasserot 
(\cite{VV09}) and Rouquier (\cite{R11}) proved that
the {\em upper global basis} introduced by the author or Lusztig's {\it 
dual canonical basis} corresponds to the isomorphism classes of 
simple $R(\beta)$-modules via the isomorphism \eqref{eq:gmod}.

However, for a given generalized Cartan matrix, 
associated \KLRs\ are not unique and depend on the parameters $c$.
 Varagnolo-Vasserot and Rouquier have proved the above results
for a very special choice $c_0$ of parameters (see \eqref{eq:Q}).
Let $R(\beta)_{c_0}$ denote the \KLR\ with the choice $c_0$, 
and $R(\beta)_\cg$ the 
\KLR\ with a generic choice $\cg$ of parameters.
When a simple $R(\beta)_\cg$-module 
is specialized at the special parameter $c_0$,
it may be a reducible $R(\beta)_{c_0}$-module.
The purpose of this note is to prove that
the specialization of any simple $R(\beta)_\cg$-module 
at $c_0$ remains a simple $R(\beta)_{c_0}$-module.
In other words, 
the set of isomorphism classes of simple $R(\beta)_\cg$-modules also 
corresponds to the upper global basis.



\noi
{\it Acknowledgements.} 
We thank 
Shunsuke Tsuchioka for helpful discussions.


\section{Review on global bases and \KLRs} \label{sec:R}
\subsection{Global bases}
Let $I$ be a finite index set. An integral square matrix
$A=(a_{i,j})_{i,j \in I}$ is called a {\em symmetrizable generalized
Cartan matrix} if it satisfies (i) $a_{i,i} = 2$ $(i \in I)$, (ii)
$a_{i,j} \le 0$ $(i \neq j)$, (iii) $a_{i,j}=0$ if $a_{j,i}=0$ $(i,j \in I)$,
(iv) there is a diagonal matrix
$D=\text{diag} (d_i \in \Z_{> 0} \mid i \in I)$ such that $DA$ is
symmetric.

A \emph{Cartan datum} $(A,P, \Pi,P^{\vee},\Pi^{\vee})$ consists of
\begin{enumerate}[(1)]
\item a symmetrizable generalized Cartan matrix $A$,
\item a free abelian group $P$ of finite rank, called the \emph{weight lattice},
\item $P^{\vee}\seteq\Hom(P, \Z)$, called the \emph{co-weight lattice},
\item $\Pi= \set{\alpha_i }{i \in I}\subset P$, called
the set of \emph{simple roots},
\item $\Pi^{\vee}= \set{h_i}{i \in I}\subset P^{\vee}$, called
the set of \emph{simple coroots},
\end{enumerate}
satisfying the condition: 
$\langle h_i,\alpha_j \rangle = a_{i,j}$ for all $i,j \in I$.

Since ${A}$ is symmetrizable,
there is a symmetric bilinear form $( \ \mid \ )$ on $P$ satisfying
$$(\alpha_i | \alpha_j)= d_i a_{i,j} \quad \text{ and } \quad  (\alpha_i | \lambda) = d_i \langle h_i, \lambda
\rangle \quad \text{ for all } i,j \in I, \ \lambda \in P.$$

The free abelian group $\rootl= \soplus_{i \in I} \Z \alpha_i$ is called the
\emph{root lattice}. Set $\rootl^{+}= \sum_{i \in I} \Z_{\ge 0}
\alpha_i\subset\rootl$ and $\rootl^{-}= \sum_{i \in I} \Z_{\le0}
\alpha_i\subset\rootl$. For $\beta=\sum_{i\in I}m_i\al_i\in\rootl$,
we set
$\haut(\beta)=\sum_{i\in I}|m_i|$.

Let $q$ be an indeterminate. Set $q_i = q^{d_i}$ for $i\in I$ and 
we define
$[n]_i =(q^n_{i} - q^{-n}_{i})(q_{i} - q^{-1}_{i} )^{-1}$ and
$[n]_i! = \prod^{n}_{k=1} [k]_i$ for $n\in\Z_{\ge0}$.

\Def\label{Def: GKM}
The {\em quantum algebra} $U_q(\g)$ associated
with a Cartan datum $({A},{P},\Pi,\Pi^{\vee})$ is 
the algebra over $\Q(q)$ generated by $e_i,f_i$ $(i \in I)$ and
$q^{h}$ $(h \in {P}^{\vee})$ satisfying following relations:
\bnum
  \item  $q^0=1$, $q^{h} q^{h'}=q^{h+h'} $ for $ h,h' \in {P}^{\vee},$
  \item  $q^{h}e_i q^{-h}= q^{\langle h, \alpha_i \rangle} e_i$,
         $q^{h}f_i q^{-h} = q^{-\langle h, \alpha_i \rangle }f_i$ 
for $h \in {P}^{\vee}, i \in I$,
  \item  $e_if_j - f_je_i =  \delta_{i,j} \dfrac{K_i -K^{-1}_i}{q_i- q^{-1}_i }, \ \ \mbox{ where } K_i=q_i^{ h_i},$
  \item  $\displaystyle \sum^{1-a_{i,j}}_{r=0} (-1)^re^{(1-a_{i,j}-r)}_i
         e_j e^{(r)}_i =0 \quad \text{ if $i \ne j$,} $
where $e_i^{(n)}=e_i^n/[n]_i!$,
  \item $\displaystyle \sum^{1-a_{i,j}}_{r=0} 
(-1)^rf^{(1-a_{i,j}-r)}_if_jf^{(r)}_i=0 \quad \text{ if $i\not=j$,}$
where $f_i^{(n)}=f_i^n/[n]_i!$.
\end{enumerate}
\enDef

Let $U_q^{-}(\g)$ be the $\Q(q)$-subalgebra of $U_q(\g)$ generated by the elements
$f_i$. We define the endomorphisms $e_i'$ and $e_i''$ of $U_q^{-}(\g)$
by
$$[e_i,a]=(q_i-q_i^{-1})^{-1}(K_ie_i''a-K_i^{-1}e_i'a)
\quad\text{for $a\in U_q^{-}(\g)$.}$$
Then $e_i'$ and the left multiplication of $f_j$
satisfy the $q$-boson commutation relations 
$$e'_if_j-q_i^{-a_{i,j}}f_je_i'=\delta_{i,j}.$$

Set $\A=\Z[q,q^{-1}]$ and let $U_\A^-(\g)$ be the $\A$-subalgebra
of $U_q^{-}(\g)$ generated by the elements
$f_i^{(n)}$.
Then $U_\A^-(\g)$ has a weight decomposition
$U_\A^-(\g)=\soplus_{\beta\in \rootl^-}U_\A^-(\g)_{\beta}$
where $U_\A^-(\g)_{\beta}\seteq\set{a\in U_\A^-(\g)}{q^haq^{-h}=q^{\lan h_i,\beta\ran}a}$.
Set $U_\A^-(\g)^*=\soplus_{\beta\in \rootl^-}\Hom_\A(U_\A^-(\g)_{\beta},\A)$
and let $e_i$, $f_i'\in\End_\A(U_\A^-(\g)^*)$ be the transposes of
$f_i, e_i'\in\End_\A(U_\A^-(\g))$, respectively.
Note that $U_\A^-(\g)_0$ is a free $\A$-module with a basis $1$,
and hence $U_\A^-(\g)^*_0$ is a free $\A$-module 
generated by the dual basis
of $1$, which is denoted by $\vac$. 

\Prop[\cite{Kash91,Kash93}]\label{Prop:glbal}
There exists a unique basis $\{\G(b)\}_{b\in B}$ of
the $\A$-module $U_\A^-(\g)^*$, called {\em the upper global basis}, which 
satisfies the following conditions:
\bnum
\item $\vac\in\set{\G(b)}{b\in B}$,
\item for any $b\in B$,
$G(b)$ belongs to $\bl U_\A^-(\g)_\beta\br^*$ 
for some $\beta\in\rootl^-$, which is denoted by $\wt(b)$,
\item Set $\eps_i(b)=\max\set{n\in\Z_{\ge0}}{e_i^n\G(b)\not=0}$.
Then for any $b\in B$ and $i\in I$, there exists $\tf_ib\in B$ such that,
 when writing 
$$f_i'\G(b)=\sum_{b'\in B}F^i_{b,b'}\G(b')\quad \text{with $F^i_{b,b'}\in \A$,}$$
we have
\bna
\item $F^i_{b,\tf_ib}=q_i^{-\eps_i(b)}$,\label{prt1}
\item $\eps_i(\tf_ib)=\eps_i(b)+1$,
\item $F^i_{b,b'}=0$ if $b'\not=\tf_ib$ and $\eps_i(b')\ge\eps_i(b)+1$,
\item $F^i_{b,b'}\in qq_i^{-\eps_i(b)}\Z[q]$ for $b'\not=\tf_ib$.\label{prt2}
\ee
\item
for $b\in B$ such that $\eps_i(b)>0$,
there exists $\te_ib\in B$ such that, when writing
$$e_i\G(b)=\sum_{b'\in B}E^i_{b,b'}\G(b')\quad \text{with $E^i_{b,b'}\in \A$,}
$$
we have
\bna
\item $E^i_{b,\te_ib}=[\eps_i(b)]_i$,
\item $\eps_i(\te_ib)=\eps_i(b)-1$, 
\item $E^i_{b,b'}=0$ if $b'\not=\te_ib$ and $\eps_i(b')\ge\eps_i(b)-1$,
\item any $E_{b,b'}^i$ is invariant under the automorphism
$q\mapsto q^{-1}$,
\item $E^i_{b,b'}\in qq_i^{1-\eps_i(b)}\Z[q]$ for $b'\not=\te_ib$.
\ee
\item $\tf_i\te_ib=b$ if $\eps_i(b)>0$,
and $\te_i\tf_ib=b$.
\ee
\enprop
Note that $B$ has the weight decomposition
$$B=\bigsqcup\nolimits_{\beta\in\rootl^-}B_{\beta}\quad\text{
with $B_{\beta}\seteq\set{b\in B}{\wt(b)=\beta}$.}
$$
There exists a unique involution (called the {\em bar involution})
$-\cl U_\A^-(\g)^*\to U_\A^-(\g)^*$
such that
\eq&&
\hs{-30ex}\parbox{40ex}{\bna
\item
$(qu)^-=q^{-1}\ol{u}$\quad for any $u\in U_\A^-(\g)^*$,
\item $-\circ e_i=e_i\circ -$\quad for any $i$,
\item $\ol{\vac}=\vac$.
\ee}
\label{char:bar}
\eneq
We have
$$\ol{\G(b)}=\G(b)\quad\text{for any $b\in B$.}$$

\subsection{Quiver Hecke algebras} 

Let $(A, P, \Pi, P^{\vee}, \Pi^{\vee})$ be a Cartan datum. In this
subsection, we recall the construction of the \KLRs\ 
associated with $(A, P, \Pi, P^{\vee}, \Pi^{\vee})$.
For $i,j\in I$ such that $i\not=j$, set
$$S_{i,j}=\set{(p,q)\in\Z_{\ge0}^2}{(\al_i,\al_i)p+(\al_j,\al_j)q=-2(\al_i,\al_j)}.
$$
Let $\cora$ be the commutative $\Z$-algebra
generated by indeterminates $\{t_{i,j;p,q}\}$
and the inverse of $t_{i,j;-a_{i,j},0}$ where 
$i,j\in I$ such that $i\not=j$ and $(p,q)\in S_{i,j}$.
They are subject to the defining relations:
$$t_{i,j;p,q} = t_{j,i;q,p}.$$

Let us define the polynomials $(Q_{ij})_{i,j\in I}$ in $\cora[u,v]$
by
\begin{equation}
Q_{ij}(u,v) = \begin{cases}\hs{5ex} 0 \ \ & \text{if $i=j$,} \\
\sum\limits_{(p,q)\in S_{i,j}}
t_{i,j;p,q} u^p v^q\quad& \text{if $i \neq j$.}
\end{cases}
\end{equation}
They satisfy $Q_{i,j}(u,v)=Q_{j,i}(v,u)$.

We denote by
$S_{n} = \langle s_1, \ldots, s_{n-1} \rangle$ the symmetric group
on $n$ letters, where $s_i\seteq (i, i+1)$ is the transposition of $i$ and $i+1$.
Then $S_n$ acts on $I^n$.

\Def[\cite{KL09,{R08}}] \label{def:KLRalg}
The {\em \KLR}\ $R(n)$ of degree $n$
associated with a Cartan datum $(A, P, \Pi, P^{\vee}, \Pi^{\vee})$ 
is the $\Z$-graded algebra over $\cora$
generated by $e(\nu)$ $(\nu \in I^{n})$, $x_k$ $(1 \le k \le n)$,
$\tau_l$ $(1 \le l \le n-1)$ satisfying the following defining
relations:
{\allowdisplaybreaks
\begin{align*}
& e(\nu) e(\nu') = \delta_{\nu, \nu'} e(\nu), \ \
\sum\nolimits_{\nu \in I^{n}}  e(\nu) = 1, \\
& x_{k} x_{l} = x_{l} x_{k}, \ \ x_{k} e(\nu) = e(\nu) x_{k}, \\
& \tau_{l} e(\nu) = e(s_{l}(\nu)) \tau_{l}, \ \ \tau_{k} \tau_{l} =
\tau_{l} \tau_{k} \ \ \text{if} \ |k-l|>1, \\
& \tau_{k}^2 e(\nu) = Q_{\nu_{k}, \nu_{k+1}} (x_{k}, x_{k+1})
e(\nu), \\
& (\tau_{k} x_{l} - x_{s_k(l)} \tau_{k}) e(\nu) = \begin{cases}
-e(\nu) \ \ & \text{if} \ l=k, \nu_{k} = \nu_{k+1}, \\
e(\nu) \ \ & \text{if} \ l=k+1, \nu_{k}=\nu_{k+1}, \\
0 \ \ & \text{otherwise},
\end{cases} \displaybreak[3]\\[.5ex]
& (\tau_{k+1} \tau_{k} \tau_{k+1}-\tau_{k} \tau_{k+1} \tau_{k}) e(\nu)\displaybreak[0]\\
&\hs{8ex} =\begin{cases} \dfrac{Q_{\nu_{k}, \nu_{k+1}}(x_{k},
x_{k+1}) - Q_{\nu_{k}, \nu_{k+1}}(x_{k+2}, x_{k+1})} {x_{k} -
x_{k+2}}e(\nu) \ \ & \text{if} \
\nu_{k} = \nu_{k+2}, \\
0 \ \ & \text{otherwise}.
\end{cases}
\end{align*}
}
The $\Z$-grading on $R(n)$ is given by
\begin{equation} \label{eq:Z-grading}
\deg e(\nu) =0, \quad \deg\; x_{k} e(\nu) = (\alpha_{\nu_k}
| \alpha_{\nu_k}), \quad\deg\; \tau_{l} e(\nu) = -
(\alpha_{\nu_l} | \alpha_{\nu_{l+1}}).
\end{equation}
\enDef

\noindent
Note that $R(n)$ has an anti-involution $\psi$ that fixes the
generators $x_k$, $\tau_l$ and $e(\nu)$.

%

For $n\in \Z_{\ge 0}$ and $\beta \in \rootl^{+}$ such that $\haut(\beta)=n$, we set
$$I^{\beta} = \set{ \nu = (\nu_1, \ldots, \nu_n) \in I^n }%
{\alpha_{\nu_1} + \cdots + \alpha_{\nu_n} = \beta }.$$
We define
\eq
&&\ba{l}
e(\beta) = \sum_{\nu \in I^{\beta}} e(\nu), \\[1ex]
 R(\beta) = R(n) e(\beta)=\soplus_{\nu\in I^\beta}R(n)e(\nu).
\ea\eneq

The algebra $R(\beta)$ is called the {\it {\KLR} at $\beta$}.

Similarly, for $\beta,\gamma\in \rootl^+$ 
with $m=\haut(\beta)$ and $n=\haut(\gamma)$
\eq
&&\ba{l}
e(\beta, \gamma) = \sum_{\nu}e(\nu)\in R(m+n)\\[1ex]
\hs{10ex}\parbox{50ex}{where $\nu$ ranges over the set of $\nu\in I^{m+n}$
such that
$\sum_{k=1}^m\alpha_{\nu_k}=\beta$ and 
$\sum_{k=m+1}^{m+n}\al_{\nu_k}=\gamma$.}
\ea\eneq
Then $R(m+n)e(\beta,\gamma)$ is 
a graded $(R(\beta+\gamma), R(\beta)\otimes R(\gamma))$-bimodule.
For a graded $R(\beta)$-module $M$ and a graded $R(\gamma)$-module $N$, 
we define their convolution
$M\circ N$ by
$$M\circ N=R(\beta+\gamma)e(\beta,\gamma)
\tens_{R(\beta)\otimes R(\gamma)}(M\otimes N). $$

For $\ell\in\Z_{\ge0}$,
We define the graded $R(\ell\al_i)$-module $\Ln(i^\ell)$ 
by 
$$\Ln(i^\ell)=q_i^{\ell(\ell-1)/2}\Bigl(
R(\ell\al_i)/\bigl(\ssum_{k=1}^\ell R(\ell\al_i)x_k\bigr)\Bigr).$$
Here $q\cl \Mod\bl R(\beta)\br\to \Mod\bl R(\beta)\br$ is the grade-shift functor:
\eq (qM)_k=M_{k-1},\eneq
and $q_i=q^{(\al_i\vert\al_i)/2}$.

For a commutative ring $\cor$ and a ring homomorphism
$c\cl \cora \to \cor$,
we denote by
$R(\beta)_\cor$ the algebra
$\cor\otimes_{\cora}R(\beta)$. 

Let us denote by $\Par$ the scheme $\mathrm{Spec}(\cora)$. For $x\in \Par$, let us denote by
$\cor(x)$ the residue field of the local ring
$(\sho_\Par)_x$
and denote by $R(\beta)_x$ the $\cor(x)$-algebra $\cor(x)\otimes_\cora R(\beta)$.

Let us take a commutative field $\cor$ and a homomorphism $\cor(x)\to\cor$.
For $\beta\in\rootl^+$, let us denote by $R(\beta)_\cor\gmod$ 
the abelian category of graded
$R(\beta)_\cor$-modules finite-dimensional over $\cor$.
Then the set of isomorphism classes of simple objects in 
$R(\beta)_\cor\gmod$ is isomorphic to the one for
$R(\beta)_{x}\gmod$ by $S\mapsto \cor\tens_{\cor(x)}S$ 
(see \cite[Corollary 3.19]{KL09}).

For $i\in I$ and $x\in \Par$ we have functors
\eqn
\xymatrix{R(\beta)_x\gmod\ar@<.5ex>[r]^-{F_i}&R(\beta+\al_i)_x\gmod
\ar@<.5ex>[l]^-{E_i}.}
\eneqn
Here these functors are defined by
\eqn
&&F_iM=
M\circ\Ln(i)\simeq\bl R(\beta+\al_i)_xe(\beta,\al_i)/R(\beta+\al_i)_xe(\beta,\al_i)x_{n+1}
\br
\otimes_{R(\beta)_x}M,\\
&&E_iN=e(\beta,\al_i)N
\simeq \Hom_{R(\beta+\al_i)_x}\bl R(\beta+\al_i)_xe(\beta,\al_i),N\br\\
&&\hs{5ex}\simeq e(\beta,\al_i)R(\beta+\al_i)_x\otimes _{R(\beta+\al_i)_x}N
\eneqn
for $M\in R(\beta)_x\gmod$ and $N\in R(\beta+\al_i)_x\gmod$.
Then we have
\eqn
&&E_i F_i\simeq q^{-(\al_i,\al_i)}F_i E_i\soplus \id,\\
&&E_i F_j\simeq q^{-(\al_i,\al_j)}F_j E_i\quad\text{for $i\not=j$,}
\eneqn
which immediately follows from \cite[Theorem 3.6]{KK}.

Let $\K{R(\beta)_x\gmod}$ denote the Grothendieck group of 
the abelian category $R(\beta)_x\gmod$.
Then, it has a structure of a $\Z[q,q^{-1}]$-module
induced by the grade-shift functor on $R(\beta)_x\gmod$.

Then the following theorem holds.
\Th[\cite{KL09}]
There exists a unique $\Z[q,q^{-1}]$-linear isomorphism
\eq
&&\soplus_{\beta\in\rootl^+}\K{R(\beta)_x\gmod}\isoto U^-_\A(\g)^*
\label{corr:main}
\eneq
such that
\bnum
\item
the induced actions $[E_i]$ and $[F_i]$ by
$E_i$ and $F_i$ correspond to
$e_i$ and $f'_i$,
\item
$\vac\in U^-_\A(\g)^*$ corresponds to the regular representation of $R(0)_x$.
\ee
\entheorem

Let $\Dual\cl R(\beta)_x\gmod\to \bl R(\beta)_x\gmod\br^\op$ be the duality functor
$M\mapsto M^*$ induced by the antiautomorphism $\psi$
of $R(\beta)_x$. 
We can easily see by the characterization \eqref{char:bar} 
of the bar involution
that the induced endomorphism
$[\Dual]$ of $\soplus\nolimits_{\beta\in\rootl^+}\K{R(\beta)_x\gmod}$ corresponds to
the bar involution $-$ of $U^-_\A(\g)^*$.

The Grothendieck group $\K{R(\beta)_x\gmod}$ is a free $\Z$-module
with the basis consisting of $[S]$ where $S$ ranges over the set of
isomorphism classes of simple graded $R(\beta)_x$-modules.
Khovanov-Lauda (\cite{KL09}) proved that
for any simple graded $R(\beta)_x$-module $S$, there exists $r\in\Z$ such that
$\Dual(q^rS)\simeq q^rS$. 
Let $\Irr(R(\beta)_x)$ be the set of isomorphism classes of
simple graded $R(\beta)_x$-modules $S$ such that $\Dual(S)\simeq S$.
Then $\K{R(\beta)_x\gmod}$ is a free $\Z[q,q^{-1}]$-module
with $\set{[S]}{S\in\Irr(R(\beta)_x)}$ as a basis.

For a simple graded module $S$, 
let us denote by $\eps_i(S)$ the largest integer
$k$ such that $E_i^kS\not=0$. 
Recall that $q$ denotes the shift-functor and $q_i=q^{(\al_i\vert\al_i)/2}$.
\Prop[\cite{LV09, KL09}]
Let $x\in \Par$, $\beta\in \rootl^+$
and  $S$ a simple graded $R(\beta)_x$-module.
\bnum 
\item
The cosocle of $F_iS$ is a simple module.
Its image under $q_i^{\eps_i(S)}$ is denoted by $\tF_iS$.
\item
If $\eps_i(S)>0$ then the socle of $E_iS$ is simple.
Its image under $q_i^{1-\eps_i(S)}$ is denoted by $\tE_iS$.
\item $\tF_i\tE_iS\simeq S$ if $\eps_i(S)>0$, and
$\tE_i\tF_iS\simeq S$.
\item If $S$ is invariant by the duality $\Dual$, then
so are $\tF_iS$ and $\tE_iS$.
\item
The set $\bigsqcup_{\beta\in\rootl^+}\Irr(R(\beta)_x)$ is isomorphic to
$B$, and 
$\tE_i$ and $\tF_i$ correspond to
$\te_i$ and $\tf_i$ by this isomorphism.
\ee
\enprop
Hence, the cosocle of $F_iS$ is isomorphic to $q_i^{-\eps_i(S)}\tF_iS$,
the socle of $E_iS$ is isomorphic to $q_i^{\eps_i(S)-1}\tE_iS$
and the cosocle of $E_iS$ is isomorphic to $q_i^{-\eps_i(S)-1}\tE_iS$.

For $b\in B_{-\beta}$, let us denote by $L_x(b)$ the corresponding simple graded
$R(\beta)_x$-module in $\Irr(R(\beta)_x)$.

Now assume that $A$ is symmetric and
consider a $\cor$-valued point
$c_0$ of $\Par$
given by
\eq
&&Q_{i,j}(u,v)=b_{i,j}(u-v)^{-a_{i,j}}\quad\text{for $i\not=j$}
\label{eq:Q}
\\
&&\hs{25ex}
\text{where $\cor$ is a field of characteristic $0$
and $b_{i,j}\in\cor^\times$.}\nonumber
\eneq
Then the following theorem is proved by Varagnolo-Vasserot
(\cite{VV09}) and Rouquier (\cite{R11}).
\Th
Assume that the generalized Cartan matrix $A$ is symmetric.
Then the basis $\{[L_{c_0}(b)]\}_{b\in B}$
corresponds to the upper global basis $\{\G(b)\}_{b\in B}$
by the isomorphism
$\soplus\nolimits_{\beta\in\rootl^+}\K{R(\beta)_{c_0}\gmod}\isoto U^-_\A(\g)^*$.
\entheorem
For $M\in R(\beta)_x\gmod$, let us define its character $\ch(M)$
by
$$\ch(M)=\sum_{\nu\in I^\beta,\;k\in\Z}\dim\bl e(\nu)M\br_k q^k e(\nu)
\in\soplus_{\nu\in I^\beta}\Z[q,q^{-1}]e(\nu).$$
Then we have
\eq
\ch\bl L_{c_0}(b)\br=\sum_{\nu\in I^\beta}\bl e_{\nu_1}\cdots e_{\nu_n}\G(b)\br
e(\nu)\quad\text{for $b\in B_{-\beta}$.}
\label{eq:character}
\eneq

\section{Main results}
\subsection{}
Let $\cg$ be the generic point of $\Par$.
For $\beta\in \rootl^+$ and $b\in B_{-\beta}$, let us consider
the simple graded $R(\beta)_{\cg}$-module $L_\cg(b)$.
\Prop
The set $U_b\seteq\set{x\in \Par}%
{\ch \bl L_x(b)\br=\ch \bl L_\cg(b)\br}$
is a Zariski open subset of $\Par$ and there exists
a graded $\sho_{U_b}\otimes_{\cora}R(\beta)$-module
$\L(b)$ defined on $U_b$ such that it is locally free as an $\sho_{U_b}$-module 
and the stalk of $\L(b)$ at any $x\in U_b$ 
is isomorphic to $L_x(b)$.
\enprop
\Proof
We shall prove it by induction on $\haut(\beta)$.
We may assume $\beta\not=0$.
Take an $i\in I$ such that $\ell\seteq\eps_i(b)\not=0$.
Set $\beta'=\beta-\ell\al_i$ and $b'=\te_i^\ell b$.
For any $x\in \Par$, 
the graded $R(\beta)_x$-module $L_x(b)$ is a simple cosocle of 
$L_x(b')\circ \Ln(i^\ell)$.
Moreover the kernel of
$L_x(b')\circ \Ln(i^\ell)\epi L_x(b)$ is 
$\set{s\in L_x(b')\circ \Ln(i^\ell)}{e(\beta',\ell\al_i)
R(\beta)s=0}$.

By the induction hypothesis,
there exists an $\sho_{U_{b'}}\otimes_{\cora}R(\beta')$-module
$\L(b')$ as above. 
Set $\R=\sho_{U_{b'}}\otimes_{\cora}R(\beta)$ and we
 shall denote by $\M$ the
$\R$-module $\L(b')\circ \Ln(i^\ell)$.
Let $f$ be the composition
\eqn
\M&\To& 
\hhom[{\sho_{\Par}\vert_{U_{b'}}}]
(\R,\M)\\
&\To&\hhom[{\sho_{\Par}\vert_{U_{b'}}}]
(\R,\M/(1-e(\beta',\ell\al_i))\M)).
\eneqn
Then the kernel of $f$ coincides with the sheaf
$$\set{u\in \M}{e(\beta',\ell\al_i)\R u=0}.$$
The homomorphismt $f$ factors through
\eqn
\M&\To[\ol{f}] &
\hhom[{\sho_{U_{b'}}}]
(\R/\R_{\ge m},\M/(1-e(\beta',\ell\al_i))\M)\br\\
&\mono&\hhom[{\sho_{U_{b'}}}]
(\R,\M/(1-e(\beta',\ell\al_i))\M)\br
\eneqn
for a sufficient large integer $m$.
Here $\R_{\ge m}=\soplus_{k\ge m}\R_k$.
Therefore $\ol{f}$ is a morphism of vector bundles on $U_{b'}$.
On the other hand, $U_b$ is the set of $x\in U_{b'}$ such that
the rank of $\ol{f}$ at $x$ is equal to its rank at the generic point.
Hence $U_b$ is an open subset of $\Par$
and the image of  $\bar{f}\vert_{U_b}$ satisfies the condition for $\L(b)$.
\QED

\subsection{}
For $x\in \Par$ and $b\in B$, let us consider the condition
\eq
&&\text{%
$L_{x}(b)$ corresponds to the upper global basis $\G(b)$ 
by the isomorphism \eqref{corr:main}.}
\label{cond:glob}
\eneq
In this subsection, we shall prove the following theorem.

\Th\label{th:main}
Let $c_0$ be a point of $\Par$ satisfying \eqref{cond:glob}
for any $b\in B$.
Then $c_0$ belongs to $U_b$ for any $b\in B$.
Hence \eqref{cond:glob} holds also for any $x\in U_b$.
\entheorem
\Proof
It is enough to show that $\cg$ satisfies  \eqref{cond:glob}.
W shall take 
a triple $(K,  \sho,\cor)$
such that $K=\cor(\cg)$, 
$\CO$ is a discrete valuation ring,
$K$ coincides with the fraction field of $\CO$, 
$\cor$ is the residue field of $\CO$,
$(\sho_\Par)_{x_0}\subset\CO$ and 
$(\sho_\Par)_{x_0}\subset\CO\to\cor$ factors through $\cor(x_0)$.
Such a triple exists (see \cite[(7.1.7)]{EGA4}).

We have the reduction map
$$\Res_{K,\cor}\cl \K{R(\beta)_K}\To \K{R(\beta)_\cor}$$
by assigning $[K\otimes_\CO L] \in \K{R(\beta)_K}$
to $[\cor\otimes_\sho L]\in\K{\R(\beta)_\cor}$
for a graded $R(\beta)_\CO$-module $L$ 
that is finitely generated and torsion-free as an $\CO$-module.
The homomorphism $\Res_{K,\cor}$ commutes with the duality
$\Dual$. Also it is compatible with the correspondence \eqref{corr:main},
namely we have a commutative diagram:
$$\xymatrix@R=2ex@C=6ex{
\soplus_{\beta\in\rootl^+}\K{R(\beta)_K\gmod}
\ar[rr]^{\Res_{K,\cor}}\ar[dr]^-\sim&&
\soplus_{\beta\in\rootl^+}\K{R(\beta)_\cor\gmod}\ar[dl]^-\sim\\
&U^-_\A(\g)^*}
$$
For $b\in B$, set $L(b)_K\seteq L_\cg(b)$ and 
$L(b)_\cor\seteq \cor\otimes_{\cor(c_0)}L_{c_0}(b)$.
Take $b\in B_{-\beta}$, and let $L(b)_\CO$ 
be an $R(\beta)_\CO$-lattice of $L(b)_K$,
i.e., a finitely generated graded $R(\beta)_\CO$-submodule $L(\beta)_\CO$ of $L(b)_K$
such that $K\otimes_\CO L(b)_\CO=L(b)_K$.
In order to see the theorem, it is enough to show that
$\cor\otimes_\CO L(b)_\CO\simeq L(b)_\cor$.

We shall prove it by induction on $\haut(\beta)$.
Take an $i\in I$ such that $\eps_i(b)>0$ and
set $b'=\te_i b$.
Then $[L(b')_K]$ corresponds to $\G(b')$ by the induction hypothesis.
We take an $R(\beta')_\CO$-lattice $L(b')_\CO$ of $L(b')_K$.
Then by the induction hypothesis, we have
$L(b')_\cor\simeq\cor\otimes_{\CO}L(b')_\CO$.
The image of $q_i^{\eps_i(b')}F_iL(b')_\CO$ by $q_i^{\eps_i(b')}F_iL(b')_K\epi L(b)_K$
is an $R(\beta)_\CO$-lattice of $L(b)_K$,
and we can take it as $L(b)_\CO$.
Since $q_i^{\eps_i(b')}F_iL(b')_\cor\simeq q_i^{\eps_i(b')}\cor\otimes_\CO 
F_i L(b')_\CO\epi\cor\otimes_\CO 
L(b)_\CO$,
the simples in a Jordan-Holder series
of $\cor\otimes_\CO L(b)_\CO$ appears in the one of $q_i^{\eps_i(b')}F_iL(b')_\cor$.

Now assume that 
$q^rL(b_1)_\cor$ appears in $\Res_{K,\cor}L(b)_K
=[\cor\otimes_{\CO}L(b)_\CO]$
for $r\in \Z$ and
$b_1\in B_{-\beta}$. 
Then $q^r\G(b_1)$ appears in 
$q_i^{\eps_i(b')}f_i'\G(b')$ by the assumption that $c_0$ satisfies 
\eqref{cond:glob}.
In particular, $L(b)_\cor$ appears in $[\cor\otimes_{\CO}L(b)_\CO]$
exactly once by  \eqref{prt1} in Proposition~\ref{Prop:glbal}.
Now assume that $(r,b_1)\not=(0,b)$.
Then \eqref{prt1} and \eqref{prt2} 
in Proposition~\ref{Prop:glbal} imply that
$r>0$. Since $L(b)_K$ is stable by the duality functor
$\Dual$, $q^{-r}L(b_1)_\cor\simeq\Dual \bl q^rL(b_1)_\cor\br$ also
appears in $\Res_{K,\cor}L(b)_K$. Hence $-r>0$.
It is a contradiction. This shows the desired result:
$\cor\otimes_\CO L(b)_\CO\simeq L(b)_\cor$.
This completes the proof of Theorem~\ref{th:main}.
\QED

\Ex
Let us give an example of a simple $R(\beta)$-module 
which does not correspond to any element in the upper global basis.
Let $\g=A^{(1)}_1$ with $I=\{0,1\}$,
$(\al_0\vert\al_0)=(\al_1\vert\al_1)=-(\al_0\vert\al_1)=2$,
and $Q_{0,1}(u,v)=u^2+auv+v^2$. Here $\cor$ is an arbitrary field
and $a\in\cor$. Set $\delta=\al_0+\al_1$,
$b'=\tf_1\tf_0\vac$
and $N=L(b')$. Then $N=\cor v$ with $x_1v=x_2v=\tau_1v=0$ and $v=e(01)v$.
Set $M=N\circ N$, and $u=v\otimes v\in M$.
Then 
$\ch(M)=2e(0101)+[2]^2e(0011)$. 
Here $e(0101)M=\cor u\oplus\cor w$ with $w\seteq\tau_2\tau_3\tau_1\tau_2u$.
By the weight consideration, $\tau_ke(0101)M=0$ for $k=1,3$ 
and $x_ke(0101)M=0$ for $1\le k\le 4$.
Easy calculations show that
$\tau_2w=-a\tau_2u$. 
Hence $y\seteq w+au$ is annihilated by all $x_k$'s and $\tau_k$'s 
and $\cor y$ is an $R(2\delta)$-submodule of $M$.
Set $M_0=M/\cor y$. Then $[M_0]$ corresponds to $\G(b)$ with
$b\seteq\tf_1^2\tf_0^2\vac$.
It is easy to see that
$M_0$ is a simple $R(2\delta)$-module if $a\not=0$.
When $a=0$, $e(0011)M_0$ is a simple  $R(2\delta)$-submodule of $M_0$
and $L(b)=e(0011)M_0$. 
Note that the case \eqref{eq:Q} is when $a=\pm2$.
\enEx
\Ex
Let us give another example of a simple $R(\beta)$-module 
which does not correspond to any element in the upper global basis.
Let $\g=A^{(1)}_2$ with $I=\Z/3\Z=\{0,1,2\}$
with $(\al_i|\al_i)=2$ and $(\al_i|\al_j)=-1$ for $i\not=j$
and $Q_{i,i+1}(u,v)=a_iu+b_{i+1}v$ ($i\in I$)
with $a_i,b_i\in\cor^\times$, where $\cor$ is an arbitrary field. 
Set $\delta=\al_0+\al_1+\al_2$,
$b'=\tf_2\tf_1\tf_0\vac$
and $N=L(b')$. Then $N=\cor v$ with $x_kv=\tau_\ell v=0$ and $v=e(012)v$.
Set $M=N\circ N$ and $u=v\otimes v\in M$.
Then 
$\ch(M)=2e(012012)+[2]^3e(001122)
+[2]^2e(001212)+[2]^2e(010122)+[2]e(010212)$. 
Here $e(012012)M=\cor u\oplus\cor w$ with
$w\seteq\tau_3\tau_4\tau_5\tau_2\tau_3\tau_4\tau_1\tau_2\tau_3u$.
By the weight consideration $\tau_ke(012012)M=0$ for $k\not=3$ 
and $x_ke(012012)M=0$ for $1\le k\le 6$.
By calculations, we have
$\tau_3w=-\gamma\tau_3u$ where $\gamma=a_0a_1a_2-b_0b_1b_2$. 
Hence $y\seteq w+\gamma u$ is annihilated by all $x_k$'s and $\tau_k$'s 
and $\cor y$ is an $R(2\delta)$-submodule of $M$.
Set $M_0=M/\cor y$. Then $[M_0]$ corresponds to $\G(b)$
with $b\seteq\tf_2^2\tf_1^2\tf_0^2\vac$.
It is easy to see that
$M_0$ is a simple $R(2\delta)$-module if $\gamma\not=0$.
When $\gamma=0$, $S\seteq 
\bl 1-e(012012)\br M_0=R(2\delta)\tau_3u$ is a simple $R(2\delta)$-module
and $L(b)=S$ and $\ch(M_0/S)=e(012012)$.
Note that the case \eqref{eq:Q} corresponds to $a_0a_1a_2+b_0b_1b_2=0$.
\enEx
\Rem
If we assume
\eq
\text{the simple modules of $R(\beta)_\cg$ 
correspond to the upper global basis,}\label{cond:global}
\eneq
then $\G(b)\in \sum_{S\in \Irr(R(\beta)_x)}\Z_{\ge0}[q,q^{-1}][S]$
for any $x\in \Par$ and $b\in B$.
We can ask if this positivity assertion still holds without the assumption 
\eqref{cond:global}.
\enRem



\bibliographystyle{amsplain}



\end{document}